\documentclass[11pt]{article}
\usepackage{amssymb}
\usepackage{amsfonts}
\usepackage{amscd}
\usepackage{epsfig}
\usepackage{amsmath,amstext,amsthm,amsbsy,amssymb}

\newcommand{\bc}{{\mathbb C}}
\newcommand{\bh}{{\mathbb H}}
\newcommand{\br}{{\mathbb R}}

\newcommand{\bb}{{\mathbb B}}
\newcommand{\mtx}[4]

\begin{document}

\title{Discreteness criterion in SL(2,$\bc$) by a test map}
\author{Wensheng Cao\thanks{Supported by NSFs of China and
NSFs of Guangdong Province}\\
School of Mathematics and Computational Science, Wuyi University,\\
Jiangmen, Guangdong 529020, P.R. China\\
e-mail:  wenscao@yahoo.com.cn}

\date{}
\maketitle
\bigskip
{\bf Abstract} \quad  In the paper (Osaka J. Math. {\bf 46}:
403-409, 2009), Yang conjectured that a non-elementary subgroup $G$
of $SL(2, \bc)$ containing elliptic elements is discrete if for each
elliptic element $g\in G$ the group $\langle f , g \rangle$ is
discrete, where $f\in SL(2,\bc)$ is  a test map which is loxodromic
or elliptic.  The purpose of this paper is to give an affirmative
answer to this question.
\smallskip

\medskip
{\bf Keywords} Discrete groups, dense groups, test map, embedding

\medskip
{\bf 2000 Mathematics Subject Classifications} 30C62, 30F40, 20H10.

\newtheorem{thm}{Theorem}[section]
\newtheorem{defi}{Definition}[section]
\newtheorem{lem}{Lemma}[section]
\newtheorem{pro}{Proposition}[section]
\newtheorem{cor}{Corollary}[section]
\newtheorem{rem}{Remark}[section]
\newtheorem{conj}{Conjecture}[section]

\section{Introduction}  The discreteness of M\"obius groups is a
fundamental problem, which has been discussed by many authors. In
1976, J{\o}rgensen established the following discreteness criterion
by using the well-known J{\o}rgensen's  inequality \cite{jor}.

\medskip

\noindent{\bf Theorem J.} \quad  {\it A non-elementary subgroup $G$
of M\"obius transformations acting on $\hat{\bc}$ is discrete if and
only if for each pair of elements $f, g\in G$,  the group
$\left\langle f, g \right\rangle$ is discrete.}

This  result  shows that the discreteness of a non-elementary
M\"obius group depends on the information of all its rank two
subgroups.  The above result has been generalized by many authors by
using  information of partial rank two subgroups.  For example,
Gilman\cite{gil} and Isochenko \cite{iso}  used  each pair of
loxodromic elements, Tukia and Wang \cite{tuwa} used each pair of
elliptic elements.

 Sullivan \cite{sul} showed  that a
non-elementary and non-discrete subgroup is either dense in
$SL(2,\bc)$ or conjugate to a dense subgroup of $SL(2,\br)$. This
result gives an approach to studying the discreteness of M\"obius
groups from the topological aspect.  Mainly using Sullivan's result,
Yang \cite{yan1} obtained some generalizations  by the information
of the remaining  four kinds of rank two subgroups.

Recently, Chen \cite{chm} proposed to use a fixed M\"obius
transformation as a test map to test the discreteness of a given
M\"obius group. His result suggests that the discreteness is not a
totally interior affair of the involved group and  provides a new
point of view to the discreteness problem.    Yang \cite{yan2}
generalized some results by test maps (see Theorems 2.4-2.7) and
proposed the following conjecture.
\medskip

\begin{conj}\label{conje}
 Let $G$ be a non-elementary subgroup of
$SL(2, \bc)$ containing elliptic elements and $f$ a loxodromic
(resp. an elliptic) transformation. If for each elliptic element
$g\in G$ the group $\langle f , g \rangle$ is discrete, then $G$ is
discrete.
\end{conj}

 Yang
proved the above conjecture for the following two special cases (
Theorems 2.9, 2.11 in \cite{yan2}).

\medskip
\noindent{\bf Theorem Y1.} \quad {\it Let $G$ be a non-elementary
subgroup of $SL(2, \br)$ containing elliptic elements and $f$ a
loxodromic (resp. an elliptic) transformation. If for each elliptic
element $g\in G$ the group $\langle f , g \rangle$ is discrete, then
$G$ is discrete.}

\noindent{\bf Theorem Y2.} \quad
 {\it Let $G$ be a non-elementary subgroup of
$SL(2, \bc)$ containing elliptic elements and $f$ a loxodromic
(resp. an elliptic) transformation with $|tr^2(f)-4|<1$. If for each
elliptic element $g\in G$ the group $\langle f , g \rangle$ is
discrete, then $G$ is discrete.}

\medskip

In $SL(2,\br)$, since the trace  is real, one can find a sequence
$\{g_n\}$ of distinct elliptic elements in $G$ such that $g_n\to I$.
In fact, this is a special case ( i.e. $dim M(G)=2$)
 of \cite[Corollary 4.5.3]{chen}.  Yang mainly used this fact to
 prove Theorem Y1.

 While in $SL(2,\bc)$, Greenberg \cite{gree} gave an example such
that  $G$ is a loxodromic group and is not discrete with $dim
M(G)=3$.  This example indicates that it is nontrivial to construct
a subgroup generated by $f$ and an elliptic element in $G$ which is
non-elementary, in which one can apply  J{\o}rgensen' inequality to
obtain a contradiction.

We mention that Theorem Y2  is true  but there is a gap in the proof
of Theorem Y2. The author got an elliptic element $g=\left(
                                 \begin{array}{cc}
                                   a & b \\
                                   c & d\\
                                 \end{array}
                               \right)$ with $b\neq 0\neq c$ and
                               $$hgh^{-1}=\left(
                                                           \begin{array}{cc}
                                                             a+c\beta & -c\beta^2+(d-a)\beta+b \\
                                                             c & -c\beta+d \\
                                                           \end{array}
                                                         \right),\ \
                                                         \mbox{where}\ h=\left(
                                                                                         \begin{array}{cc}
                                                                                           1 & \beta \\
                                                                                           0 & 1 \\
                                                                                         \end{array}
                                                                                       \right).$$
Taking $\beta=\frac{d-a}{2c}$,  the author  mistook the second entry
of $hgh^{-1}$ for zero. In fact, $$hgh^{-1}=\left(
                   \begin{array}{cc}
                    \frac{a+d}{2} & \frac{(a+d)^2-4}{4c}\\
                            c & \frac{a+d}{2} \\
                               \end{array}
                                 \right).$$  If the sequence $\{h_n\}$ in $G$ converges to $h$,
                                 then the product $b_nc_n$ of the second and third entries of $h_ngh_n^{-1}
                                  $ converges to $\frac{(a+d)^2-4}{4}$ which is not
                                 zero for $g$ being elliptic.

We can mend the the proof of Theorem Y2 as followings.

  {\bf The proof of Theorem Y2.}\quad Suppose that
$G$ is not discrete.

 If $G$ is a dense subgroup in $SL(2,\br)$ then
as reasoning in Theorem Y1, we can get the result.

 If $G$ is dense in $SL(2,\bc)$, we can solve the
following equation
\begin{equation}\label{eq1}-cz^2+(d-a)z+b=0 \end{equation} to get a solution $\beta$.
Construct $h=\left(
                                               \begin{array}{cc}
                                                     1 & \beta \\
                                                       0 & 1 \\
                                                   \end{array}
                                                  \right)$ with this obtained $\beta$. Let $\{h_n\}$  be  a sequence in $G$ converges to
                                                  $h$ and $g_n=h_ngh_n^{-1}=\left(
   \begin{array}{cc}
     a_n & b_n \\
     c_n & d_n \\
   \end{array}
 \right)$.  Then $b_nc_n\to 0$ and  $ \langle f,g_n \rangle$ is discrete and non-elementary for large $n$.   This contradicts the
J{\o}rgensen's inequality
$$|tr(f)^2-4|+|tr[f,g_n]-2|=(1+|b_nc_n|)|r-\frac{1}{r}|^2\geq 1.$$

The proof is complete.

\medskip

Let $g=\left(
                                 \begin{array}{cc}
                                   a & b \\
                                   c & d\\
                                 \end{array}
                               \right)$ be an elliptic element and $f=\left(
                                                             \begin{array}{cc}
                                                               r & 0 \\
                                                               0 & \frac{1}{r} \\
                                                             \end{array}
                                                           \right)$ be a loxodromic or an elliptic element. Then
\begin{equation} |tr(g)^2-4|+|tr[g,f]-2|=4-(a+d)^2+|bc||r-\frac{1}{r}|^2.\end{equation}
Suppose $G$ is dense in $SL(2,\bc)$ and non-elementary. If we can
find an elliptic element $g= \left(
                           \begin{array}{cc}
                             a & b \\
                             c & d \\
                           \end{array}
                         \right)\in G$
 with $4-(a+d)^2<1$, then as in the above
proof, we can get a sequence $\{g_n\}$ of distinct elliptic elements
with $|b_nc_n|\to 0$.  This may provide us a desired contradiction
to prove Conjecture \ref{conje}. Motivated by this  observation, we
manage to get such an elliptic element under certain condition by
embedding $SL(2,\bc)$ into $U(1,1; \bh)$.

Our main theorem is

\begin{thm}\label{thm1}
Conjecture \ref{conje} is positive.\end{thm}

\section{The unitary group and embedding principle}
In this section, we will recall some facts about quaternion and the
quaternionic hyperbolic geometry.  The reader is referred to
\cite{cpw,cata,chen} for more information.

 Let $\bh$ denote the
division ring of real quaternions. Elements of $\bh$ have the form
$q=q_1+q_2{\bf i}+q_3{\bf j}+q_4{\bf k}\in \bh$ where $q_i\in \br$
and
$$
{\bf i}^2 = {\bf j}^2 = {\bf k}^2 = {\bf i}{\bf j}{\bf k} = -1.
$$
Let $\overline{q}=q_1-q_2{\bf i}-q_3{\bf j}-q_4{\bf k}$ be the {\sl
conjugate} of $q$, and
$$
|q|= \sqrt{\overline{q}q}=\sqrt{q_1^2+q_2^2+q_3^2+q_4^2}
$$
be the {\sl modulus} of $q$. We define $\Re(q)=(q+\overline{q})/2$
to be the {\sl real part} of $q$, and $\Im(q)=(q-\overline{q})/2$ to
be the {\sl imaginary part} of $q$. Also
$q^{-1}=\overline{q}|q|^{-2}$ is the {\sl inverse} of $q$. We remark
that for a complex number $c$, we have ${\bf j}c=\bar{c}{\bf j}$.

 Let
$\bh^{1,1}$ be the vector space of dimension 2 over $\bh$ with the
unitary structure defined by the Hermitian form
$$
\langle{\bf z},\,{\bf w}\rangle={\bf w}^*J{\bf z}=
\overline{w_1}z_1-\overline{w_2}z_2,
$$
where ${\bf z}$ and ${\bf w}$ are the column vectors in $\bh^{1,1}$
with entries $(z_1,z_2)$ and $(w_1,w_2)$ respectively, $\cdot^*$
denotes the conjugate transpose and $J$ is the Hermitian matrix
$$
J=\left(\begin{matrix} 1 & 0 \\ 0 & -1 \end{matrix} \right).
$$
We define a {\sl unitary transformation} $g$ to be an automorphism
$\bh^{1,1}$, that is, a linear bijection such that
\begin{equation}\label{model}\langle g({\bf z}),\,g({\bf w})\rangle=\langle{\bf
z},\,{\bf w}\rangle\end{equation} for all ${\bf z}$ and ${\bf w}$ in
$\bh^{1,1}$. We denote the group of all unitary transformations by
${\rm U}(1,1;\bh)$.

Following \cite[Section 2]{chen}, let
$$
V_0 =  \Bigl\{{\bf z} \in  \bh^{1,1} -\{0\}: \langle{\bf z},\,{\bf
z}\rangle=0\Bigr\}, \  \ V_{-}   =  \Bigl\{{\bf z} \in
\bh^{1,1}:\langle{\bf z},\,{\bf z}\rangle<0\Bigr\}.
$$
It is obvious that $V_0$ and $V_{-}$ are invariant under ${\rm
U}(1,1;\bh)$. We define $V^s$ to be $V^s=V_{-}\cup  V_0$. Let
$P:V^s\to P(V^s)\subset \bh$ be the projection map defined by
$$
P\left(\begin{matrix} z_1 \\ z_2
\end{matrix}\right)=z_{1}{z_2}^{-1}.
$$
We define $\bb=P(V_-)$, the ball model of 1-dimensional quaternionic
hyperbolic space. It is easy to see that $\bb$ can be identified
with the quaternionic unit ball $\bigl\{ z \in \bh:\;|z|<1\bigr\}$.
Also the unit sphere in $\bh$ is $\partial\bb=P(V_0)$ and the center
of the ball is $0=P\left(\begin{matrix} 0 \\ 1
\end{matrix}\right)$.

If $ g=\left(\begin{matrix} a & b \\ c & d \end{matrix}\right) \in
{\rm U}(1,1; \bh) $ then, by definition, $  g$ preserves the
Hermitian form. Hence
$$
{\bf w}^*J{\bf z}=\langle{\bf z},\,{\bf w}\rangle= \langle g{\bf
z},\,  g{\bf w}\rangle ={\bf w}^*  g^*J  g{\bf z}
$$
for all ${\bf z}$ and ${\bf w}$ in $V$. Letting ${\bf z}$ and ${\bf
w}$ vary over a basis for $V$ we see that $J=  g^*J  g$. From this
we find $  g^{-1}=J^{-1}  g^*J$. That is:
$$
\left(\begin{matrix} a & b \\ c & d \end{matrix}\right)^{-1}
=\left(\begin{matrix} \overline{a} & -\overline{c} \\
-\overline{b} & \overline{d} \end{matrix}\right)
$$
 and  consequently,

\begin{equation}\label{relation}
|a|=|d|, \; |b|=|c|,\; |a|^2-|c|^2=1, \;  \bar{a}b=\bar{c}d, \;
a\bar{c}=b\bar{d}.
\end{equation}

As in \cite{cpw,cata}, we can regard $U(1,1;\bh)$ as the isometries
of real hyperbolic 4-space, whose model is the unit ball in the
quaternions $\bh$.  $SL(2,\bc)$, the isometries of real hyperbolic
3-space, can be embedded as a subgroup of $U(1,1;\bh)$ as following:

$$f\in SL(2,\bc)\hookrightarrow TfT^{-1}\in U(1,1;\bh),$$
where $$T=\frac{1}{\sqrt{2}}\left(
            \begin{array}{cc}
              1 & -{\bf j} \\
              -{\bf j} & 1 \\
            \end{array}
          \right).$$

Let $f=\left(
            \begin{array}{cc}
              a & b \\
              c & d \\
             \end{array}
          \right)\in SL(2,\bc)$. Then
$$ \hat{f}=TfT^{-1}=\frac{1}{2}\left(
                                                            \begin{array}{cc}
                                                              1 & -{\bf j} \\
                                                             -{\bf j} & 1 \\
                                                            \end{array}
                                                          \right)\left(
                                                                   \begin{array}{cc}
                                                                     a & b \\
                                                                     c & d \\
                                                                   \end{array}
                                                                 \right)
                                                          \left(
                                                            \begin{array}{cc}
                                                              1 & {\bf j} \\
                                                              {\bf j} & 1 \\
                                                            \end{array}
                                                          \right)\in U(1,1;\bh).$$

We mention that our model is slight different from the model in
\cite{chen}, where  the Hermitian matrix is $ J=\left(\begin{matrix}
-1 & 0 \\ 0 & 1 \end{matrix} \right).$   It follows from
(\ref{model}) that both models define the same unitary group. This
difference just exchanges the inner and outer of the same unit
sphere of those two models.

\smallskip

The following lemma is crucial to us.

\begin{lem}\label{lem1}(cf. \cite[Corollary 4.5.2]{chen}) Let G be a subgroup of $U(1,n;\bh)$ such that (a) $G$
does not leave invariant a point in $\partial H_{\bh}^n$ or a proper
totally geodesic submanifold of $H_{\bh}^n$ (b) the identity is not
an accumulation point of the elliptic elements in $G$. Then $G$ is
discrete.
\end{lem}

Using the same notation as in \cite{chen}, for any totally geodesic
submanifold $M\in H_{\bh}^n$, we denote by $I(M)$ the subgroup of
$U(1,n;\bh)$ which leaves $M$ invariant. By \cite[Proposition
2.5.1]{chen}, the proper totally geodesic submanifolds of
$H_{\bh}^1$ are equivalent to one of the four types:  $H_{\br}^1$,
$H_{\bc}^1$ and $H^1(\mathbb{I})$.

\medskip

By \cite[Lemmas 4.2.1,2]{chen}, we have the following lemma.

\begin{lem}\label{lem2} Let $g\in U(1,1;\bh)$.   Then

(i) the elements $g\in I(H_{\br}^1)$ are of the form
$$g=A\lambda, \ A\in U(1,1;\br), \lambda\in \bh, |\lambda|=1;$$

(ii)the elements $g\in I(H_{\bc}^1)$ are of the form
$$g=A,  A\in U(1,1;\bc);$$

(iii) the elements $g\in I(H^1(\mathbb{I}))$ are of the form
\begin{equation} \label{im} g=\left(
      \begin{array}{cc}
        a & b \\
       -\varepsilon b & \varepsilon a \\
      \end{array}
    \right)\in U(1,1;\bh),
\varepsilon=\pm 1.\end{equation}
\end{lem}

\begin{lem}\label{lem3} Let $G$ be a subgroup of $SL(2,\bc)$. Then $TGT^{-1}$ is a subgroup of $U(1,1;\bh)$.
If $g=\left(
            \begin{array}{cc}
              a & b \\
              c & d \\
             \end{array}
          \right)\in G$ and $TGT^{-1}\subset I(H^1(\mathbb{I})) $
          then either

(i) $a,d\in \br$  and $b,c\in {\bf i}\br$,or

(ii) $a,d\in {\bf i}\br$  and $b,c\in \br$.
\end{lem}

{\bf Proof.}  If $g=\left(
            \begin{array}{cc}
              a & b \\
              c & d \\
             \end{array}
          \right)\in G$ and $TGT^{-1}\subset I(H^1(\mathbb{I})),$
          then
$ TgT^{-1}$  is of form (\ref{im}). By our embedding and the fact
${\bf j}c=\bar{c}{\bf j}, \forall c\in\bc$, we can verify that the
cases $\varepsilon=1$ and $\varepsilon=-1$ correspond to cases (i)
and (ii), respectively.

\begin{lem}\label{lem4} Let $G$ be a subgroup of $U(1,1;\bh)$.
Define $tr(g)=a+d$ for $g=\left(
            \begin{array}{cc}
              a & b \\
              c & d \\
             \end{array}
          \right)\in G$. Then
$$\Re(tr(g))=\Re(tr(fgf^{-1})), \forall f\in U(1,1;\bh). $$

\end{lem}

\section{The proof of Theorem \ref{thm1}}

We also need the following lemma, which is a direct consequence of
the well-known proposition in \cite[Section 1]{sul}.

\begin{lem}  Let G be a non-elementary  subgroup of
$SL(2,\bc)$. Then either

(i) $G$ is discrete,or

(ii) $G$ is  dense in $SL(2,\bc)$, or

(iii) $G$ is  conjugate to a dense group of $SL(2,\br)$.
\end{lem}

{\bf The proof of Theorem \ref{thm1}.} \quad Suppose that $G$ is
non-elementary and not discrete.    If $G$ is conjugate to a dense
group of $SL(2,\br)$ then we can obtain the result as in Theorem Y1.

In what follows,  we assume that  $G$ is  dense in $SL(2,\bc)$.

 By our
embedding, $G_1=TGT^{-1}$ is a non-elementary and non-discrete
subgroup of $U(1,1;\bh)$.
 Let $M(G_1)$ be the smallest totally geodesic submanifold which is
invariant under $G_1$. By \cite[Lemma 4.5.1]{chen}, the limit set
$L(G_1)$ of $G_1$ belongs to $\partial M(G_1)$.  Normalize such that
$M(G_1)$ contains $0$  (abuse of notation, still denote this
normalized subgroup by $G_1$), then $M(G_1)$ is one of the four
types: $H_{\br}^1$, $H_{\bc}^1$, $H^1(\mathbb{I})$ and $H_{\bh}^1$.

Since $TGT^{-1}$ is non-elementary,  $M(G_1)\neq H_{\br}^1$. Suppose
that $M(G_1)= H_{\bc}^1$.   By Lemma \ref{lem2}, $TGT^{-1}$ is a
subgroup of $U(1,1;\bc)$. Since $PU(1,1;\bc)$ is isomorphism to
$PSL(2,\br)$,  we can get the result as in Theorem Y1 in this case.

Suppose that $M(G_1)= H_{\bh}^1$. By Lemma \ref{lem1}, the identity
is  an accumulation point of the elliptic elements in $G_1$.
Therefore we get a sequence $\{g_n\}$ of distinct elliptic elements
in $G$ such that $g_n\to I$.   By the same reasoning as in  Theorem
Y1, we can get the result.

Suppose that $M(G_1)= H^1(\mathbb{I})$. By Lemmas \ref{lem3},
\ref{lem4}, we know that the trace of $g\in G$  belongs to either
$\br$ or ${\bf i}\br$.  Let $k\in SL(2,\bc)$ be an elliptic element
with $3.1<tr^2(k)<3.9$. Since $G$ is dense in $SL(2,\bc)$, there
exist a sequence $k_n$ converges to $k$.   Therefore we can find an
elliptic element $g=\left(
             \begin{array}{cc}
               a & b \\
               c & d \\
             \end{array}
           \right)
\in SL(2,\bc)$ with $3<tr^2(g)<4$. Since $G$ is non-elementary, we
can further assume that $b\neq 0\neq c$.

Let $\beta$ be a solution to the equation (\ref{eq1}).  Construct
$h=\left(
                                               \begin{array}{cc}
                                                     1 & \beta \\
                                                       0 & 1 \\
                                                   \end{array}
                                                  \right)$ and let $\{h_n\}$  be  a sequence in $G$ converges to
                                                  $h$.  Then
 $g_n=h_ngh_n^{-1}=\left(
   \begin{array}{cc}
     a_n & b_n \\
     c_n & d_n \\
   \end{array}
 \right)$ are  elliptic elements  with $3<tr^2(g_n)=tr^2(g)<4$ and  $b_nc_n\to 0$ as  $n\to \infty$.
  Note that $ \langle f,g_n \rangle$ is discrete and non-elementary for large $n$. This  contradicts the
J{\o}rgensen's inequality
\begin{equation} |tr(g_n)^2-4|+|tr[g_n,f]-2|=4-(a_n+d_n)^2+|b_nc_n||r-\frac{1}{r}|^2\geq 1\end{equation}

  The proof is complete.

\medskip

{\bf Remark.}\quad   In \cite{yan3}, Yang asked that whether there
is a non-elementary and nondiscrete subgroup of
$Isom(H^3)=PSL(2,\bc)$ which contains elliptic such that each of
them has order 2. By the proof of Theorem \ref{thm1}, we know that
the answer to this question is {\it negative}.

\end{document}